\documentclass[12pt]{amsart}
\usepackage{amssymb}
\usepackage{amsfonts}
 \usepackage{amsfonts,amssymb}

\usepackage{mathrsfs}
\let\mathcal\mathscr

\usepackage[all,ps,cmtip]{xy}

\def\llra{\hbox to 10mm{\rightarrowfill}}

\def\lllra{\hbox to 15mm{\rightarrowfill}}

\def\phi{{\varphi}}

\def\cI{\mathcal{I}}

\def\cA{\mathcal{A}}
\def\cF{\mathcal{F}}

\def\cO{\mathcal{O}}

\def\cP{\mathcal{P}}

\def\cJ{\mathcal{J}}

\def\cQ{\mathcal{Q}}

\def\cZ{\mathcal{Z}}

\let\tilde\widetilde

\DeclareMathOperator{\codim}{codim}
\DeclareMathOperator{\Pic}{Pic}

\DeclareMathOperator{\vol}{vol}

\newtheorem{lemm}{Lemma}[section]
\newtheorem{theo}[lemm]{Theorem}
\newtheorem{coro}[lemm]{Corollary}

\newtheorem{prop}[lemm]{Proposition}

\newtheorem*{conj*}{Conjecture}

\theoremstyle{definition}

\newtheorem{rema}[lemm]{Remark}

\newtheorem{conj}[lemm]{Conjecture}

\theoremstyle{remark}
\newtheorem*{remark*}{Remark}
\newtheorem*{note*}{Note}

\begin{document}
\title{ Angehrn-Siu-Helmke's method applied to abelian varieties }
\author{Zhi Jiang}
\address{Zhi Jiang, Shanghai Center for Mathematical Sciences, China}
\email{zhijiang@fudan.edu.cn}
\date{\today}
\maketitle
\begin{abstract}
We apply  Angehrn-Siu-Helmke's method to estimate basepoint freeness thresholds of higher dimensional polarized abelian varieties.  We showed that a conjecture of Caucci  holds for very general polarized abelian varieties in the moduli spaces $\cA_{g, l}$ with only finitely many possible exceptions of polarization types $l$ in each dimension $g$. We improved the bound of basepoint freeness thresholds of any polarized ableian $4$-folds and simple abelian $5$-folds.
\end{abstract}
\maketitle
\section{Introduction}
Syzygies of abelian varieties have attracted lots of attention in recent years. Recall the following question asked by Ito for $p\geq 0$ and by Lozovanu for $p=-1$ as well (see \cite{Ito1} and \cite{Loz}).
\begin{conj}\label{P1}Let $(A, L)$ be a polarized abelian variety of dimension $g$ and $p\geq -1$ an integer. If $(L^g)>(g(p+2)))^g$ and $(L^d\cdot B)>(d(p+2))^d$ for any abelian subvariety $B$ of $A$ of dimension $0<d<g$, then $L$ satisfies Property $(N_p)$.
\end{conj}
We summarize recent progress towards this conjecture.
\subsection{$\mathbb Q$-twisted sheaves} Given a coherent sheaf $\cF$ on $A$ and a rational number $t\in \mathbb Q$, following \cite{JP}, we formally define the $\mathbb Q$-twisted sheaves $\cF\langle tL \rangle$. We say that $\cF\langle tL\rangle$ is IT$^0$ if the $i$-th cohomological support loci $V^i(\mu_b^*\cF\otimes L^{\otimes b^2t})$ is empty for each $i>0$, where $b$ is an integer such that $b^2t$. Note that this definition does not depend on the choice of $b$. Assume that $D$ is an effective $\mathbb Q$-divisor on $A$ such that $D$ is $\mathbb Q$-equivalent to $tL$, we will also write $\cF\langle D\rangle=\cF\langle tL \rangle$.

Similarly, we say that $\cF\langle tL \rangle$ is M-regular (resp. GV) if $$\codim_{\Pic^0(A)}V^i(\mu_b^*\cF\otimes L^{\otimes (b^2t)})>i$$  (resp. $\codim_{\Pic^0(A)}V^i(\mu_b^*\cF\otimes L^{\otimes (b^2t)})\geq i$) for all $i>0$.
We can similarly define the cohomology ranks of $\cF\langle tL\rangle$: $h^i(A, \cF\langle tL\rangle):=\frac{1}{b^{2\dim A}}h^i(A, \mu_b^*\cF\otimes L^{\otimes b^2t}\otimes Q)$, where $Q\in\Pic^0(A)$ is general. The main result of \cite{JP} says that the function $t\rightarrow  h^i(A, \cF\langle tL\rangle)$ is locally polynomial on a left or right neighborhood of a rational number and can be extended to a continuous function from $\mathbb R$ to $\mathbb R$. We call this function the $i$-th cohomological rank function of $\cF$ and denote it by $h^i_{\cF, L}(t)$.

Let $\cI_o$ be the ideal sheaf of the neutral element $o$ of $A$. The basepoint freeness threshold $\beta(L)$ is defined to be $$\beta(L):=\inf\{t\in \mathbb Q\mid \cI_o\;\langle tL\rangle \;is\; \mathrm{IT}^0\}.$$ By the main theorem of \cite{Hac}, $\beta(L)$ is also equal to $$\inf\{t\in \mathbb Q\mid \cI_o\;\langle tL\rangle \;is\; \mathrm{GV}\}.$$
It was observed in \cite{JP} that $\beta(L)\leq 1$ and equality holds iff $|L|$ has a basepoint. This shows that $\beta(L)$ may vary in families.  More generally, if for some rational number $t=\frac{a}{b}$,   $V^i(\mu_b^*\cI_o\otimes L^{\otimes (ab)})$ is a non-empty proper subset of $\Pic^0(A)$, $\beta(L)=t$.
But there is by far no general way to determine $\beta(L)$.
It is also not clear that whether or not $\beta(L)$ is always a rational number.

By \cite{JP} and \cite{C}, knowing the exact number of the basepoint freeness threshold $\beta(L)$ helps to understand the syzygies of $L$.
\begin{theo}\label{criteria}For $p\geq -1$, if $\beta(L)<\frac{1}{p+2}$, $L$ satisfies Property $(N_p)$.
\end{theo}
\begin{rema}Ito refined this criteria in \cite{Ito4} by showing that  if for some integer $p> 0$, $\beta(L)=\frac{1}{p+2}$ and $\cI_o\langle \frac{1}{p+2}L\rangle$ is M-regular, $L$ satisfies Property $(N_p)$.
\end{rema}
Caucci then asked the following question, which would imply Conjecture \ref{P1}.
\begin{conj}\label{P2}Let $(A, L)$ be a polarized abelian variety of dimension $g$ and $p\geq -1$ an integer. If $(L^g)>(g(p+2)))^g$ and $(L^d\cdot B)>(d(p+2))^d$ for any abelian subvariety $B$ of $A$ of dimension $0<d<g$, then $\beta(L)<\frac{1}{p+2}$.
\end{conj}
\subsection{Known results}
By Theorem \ref{criteria}, in order to solve Conjecture \ref{P1} and Conjecture \ref{P2}, it suffices to prove that  \begin{eqnarray}\label{mainpoint}&&\beta(L)\leq n(L):=\\\nonumber  \mathrm{inf}\{\frac{d}{\sqrt[d]{(L^d\cdot B)}}\mid &&  B  \;\mathrm{ \; is \;an \;abelian\; subvariety\;of\; dimension}\;  1\leq d \leq g\}.\end{eqnarray}
There are various ways to estimate $\beta(L)$. In \cite{Ito}, Ito showed that the Angehrn-Siu method (see \cite{AS, Hel, Kaw}), which was initially applied to attack Fujita's basepoint-freeness conjecture, can also be used to estimate $\beta(L)$. To be more precise, let \begin{eqnarray*}r'(L):=\mathrm{inf}\{t\in \mathbb Q\mid \mathrm{there\; exists\; a\;} \mathbb{Q}-divisor D\sim_{\mathbb Q}tL \\ \mathrm{such\; that\; o\; is\; an\; isolated\; component\; of\; Nklt(A, D)}\},\end{eqnarray*} where the non-klt locus $\mathrm{Nklt}(A, D)$ is the subscheme of $A$ defined by the multiplier ideal $\cJ(A, D)$ (see for instance \cite[Section 9 and 10]{Lar2}). In \cite{Ito}, Ito observed that $\beta(L)\leq r'(L)$, when $(A, L)$ is a polarized simple abelian 3-fold, $  \beta(L)\leq r'(L)\leq n(L)$, and Conjecture \ref{P1} and Conjecture \ref{P2} hold for any polarized abelian 3-fold.


In higher dimensions,   we  have proved in \cite{J} via generic vanishing  that $\beta(L)\leq  2n(L).$ This implies that given a polarized abelian variety $(A, L)$, if $(L^d\cdot B)>(2(p+2)d)^d$ for any abelian subvariety $B$ of dimension $1\leq d\leq g$, $L$ satisfies Property $(N_p)$.

Based on Bridgeland's stability condition of surfaces, Lahoz-Rojas \cite{LR} and Rojas \cite{R}  almost determined the cohomological rank functions $h^i_{\cI_o, L}$ for any polarized abelian surface $(A, L)$ of Picard number $1$. When $L$ is of polarization type $(1, d)$, Rojas proved that $\beta(L)=\frac{1}{\sqrt{d}}$ when $d$ is a perfect square or $\beta(L)=\frac{2y}{x-1}$ where $(x, y)$ is the minimal or the second minimal positive solution of the Pell's equation $X^2-4dY^2=1$ when $d$ is not a perfect square. In either case, we have $\beta(L)< \frac{\sqrt{d}+1}{d}<\frac{\sqrt{2}}{\sqrt{d}}=n(L)$ when $d\geq 6$.

Via a degeneration method, Ito studied in \cite{Ito3} the syzygies of general polarized abelian varieties of type $(1,\ldots,1, d)$. He proved that when $d\geq \frac{(p+2)^{g+1}-1}{p+1}$, $L$ satisfies Property $(N_p)$.

Inspired by the results of Ito and Rojas, one may believe that $\beta(L)$ should be quite close to $\frac{1}{\sqrt[g]{h^0(A, L)}}=\frac{\sqrt[g]{g!}}{\sqrt[g]{(L^g)}}$ for general polarized abelian variety $(A, L)$ of dimension $g$.
\subsection{Main results}
In this paper, we follow Ito's approach to apply the Angehrn-Siu method to estimate $\beta(L)$ for polarized abelian varieties. We improve  the known upper bounds for $\beta(L)$ for higher dimensional polarized abelian varieties.

We denote respectively by $g$ and $l$ the dimension and the polarization type of a polarized abelian variety $(A, L)$. There is a quasi-projective variety $\cA_{g, l}$ parametrizing such polarized abelian varieties.
\begin{theo}\label{thm1}
Let $(A, L)$ be a very general polarized abelian  variety in $\cA_{g, l}$. When $g=4$ or $5$, we have $\beta(L)\leq n(L)$. When $g\geq 6$,  we also have $\beta(L)\leq n(L)$, except for possibly finitely many polarization types $l$.
\end{theo}

\begin{theo}\label{thm2}
Let $(A, L)$ be  a polarized abelian $4$-fold. Assume that $(L^4)>((2+\frac{4}{\sqrt{3}})(p+2))^4$ and $(L^d\cdot B)>((p+2)d)^d$ for any abelian subvariety $B$ of dimension $1\leq d\leq 3$. Then $\beta(L)<\frac{1}{p+2}$.
\end{theo}
\begin{rema} Note that $2+\frac{4}{\sqrt{3}}\approx 4.31$. Thus Theorem \ref{thm2} is quite close to Conjecture \ref{P2} for abelian 4-folds.
\end{rema}
\begin{theo}\label{thm3}
Let $(A, L)$ be  a polarized abelian $5$-fold. Assume that $A$ is simple and $(L^5)>(8(p+2))^5$. Then $\beta(L)< \frac{1}{p+2}$.
\end{theo}

In higher dimensions, we also have   slight improvements of $\beta(L)$.  For $g\geq 6$, we define $\alpha_{g, g-2}=\sqrt[g-2]{\frac{15(g-3)!}{g-1}}$, $\alpha_{g, 2}=\sqrt{\frac{5(g-2)+1}{g-1}}$, and for $3\leq d\leq g-3$, $$\alpha_{g,d}=\sqrt[d]{\frac{\frac{5}{2}d!(g-d)}{\binom{g-1}{g-d}}}$$ and let $$\alpha_g:=\mathrm{min}\{\alpha_{g,d}\mid 2\leq d\leq g-2\}.$$
\begin{theo}\label{thm4} Let $(A, L)$ be a polarized abelian variety of dimension $g\geq 6$. Assume that $A$ is simple and $(L^g)>((2g-\alpha_g)(p+2))^g$, $\beta(L)< \frac{1}{p+2}$.
\end{theo}

Note that $\alpha_6=\alpha_{6,3}=\sqrt[3]{\frac{9}{2}}$, $\alpha_7=\alpha_{7, 3}=\sqrt[3]{4}$, $\alpha_8=\alpha_{8, 3}=\sqrt[3]{\frac{25}{7}}$.

\subsection*{Acknowledgements} The author thanks Atsushi Ito for helpful communications and in particular, Lemma 2.15 is due to him. The author also thanks Chen Jiang for helpful conversations and Federico Caucci for helpful comments.
The author is supported by NSFC for Innovative Research Groups
(No. 12121001), by the National Key Research and Development Program of China (No. 2020YFA0713200), by the Natural Science Foundation of China (No. 11871155 and No. 11731004), and by the Natural
Science Foundation of Shanghai (No. 21ZR1404500).
 \section{Preliminaries}
\subsection{Log canonical centers}

Let $X$ be a smooth projective variety and  $\Delta$ be an effective $\mathbb Q$-divisor on $X$.
 We take $\mu: Y\rightarrow X$  a log resolution of $(X, \Delta)$ and we write $$K_Y=\mu^*(K_X+\Delta)+\sum_Ea(E, X, \Delta)E,$$ where  $E$ runs through prime divisors of $\mu$ and $\mu_*(\sum_Ea(E, X, \Delta)E)=-\Delta$.

We say that $(X, \Delta)$ is log canonical at $x$ if  $a(E, X, \Delta)\geq -1$ for all   prime divisors $E$ on $Y$ such that $x\in \mu(E)$. For a prime divisor $E$ on $Y$, if $(X, \Delta)$ is log canonical at the generic point $\mu(E)$ and the discrepancy $a(E, X, \Delta)=-1$, we call $\mu(E)$ a log canonical center of $(X, \Delta)$.

 When $(X, \Delta)$ is log canonical at $x$, there are finitely many log canonical centers containing $x$ and the intersection of   two such log canonical centers is the union of certain log canonical centers containing $x$
 (see for instance \cite[Proposition 1.5]{Kaw} or \cite[Theorem 9.1]{F}). Thus there exists a unique minimal log canonical center $Z$ of $(X, \Delta)$ through $x$.

 Moreover, $Z$ is normal and has rational singularities around $x$. More precisely, locally around $x$, there exists a $\mathbb Q$-effective divisor $\Delta_Z$ such that $(Z, \Delta_Z)$ is klt  (see for instance \cite[Theorem 7.1]{FG}).
However, the singularities of $Z$ away from $x$ cannot be controlled in general. The following result due to Xiaodong Jiang from \cite[Proposition 5.1]{Jx} will be applied later.
\begin{prop}\label{Jx} Let $X$ be a smooth projective variety and $\Delta$ an effective $\mathbb Q$-divisor on $X$ with $Z$  a log canonical center of $(X, \Delta)$ \footnote{Note that  the notion of  a pure log canonical center in \cite{Jx} coincides with the notion of a log canonical center here.}.
 Let $\nu: \overline{Z}\rightarrow Z$ be the normalization. Then  there exists an effective $\mathbb Q$-divisor $\Delta_{\overline{Z}}$ on $\overline{Z}$ such that $\nu^*(K_X+\Delta)\sim_{\mathbb Q} K_{\overline{Z}}+\Delta_{\overline{Z}}$.
\end{prop}

We then have the immediate corollary.
\begin{coro}\label{intersection}
 Under the above assumption, we have $$\vol(Z, (K_X+\Delta)|_Z)=\vol(\overline{Z}, K_{\overline{Z}}+\Delta_{\overline{Z}})\geq \vol(Z', K_{Z'}),$$ where $Z'$ is any $\mathbb Q$-Gorenstein partial resolution of $\overline{Z}$.
\end{coro}
\begin{proof}
Let $\rho: Z'\rightarrow Z$ be a $\mathbb Q$-Catier partial resolution and
let $\mu: \tilde{Z}\xrightarrow{\sigma }Z'\xrightarrow{\rho} \overline{Z}$ be a log resolution of $(\overline{Z}, \Delta_{\overline{Z}})$. Let $\Delta_{\tilde{Z}}=\mu_*^{-1}(\Delta_{\overline{Z}})$. We then have
$K_{\tilde{Z}}+\Delta_{\tilde{Z}}+E_1=\mu^*(K_{\overline{Z}}+\Delta_{\overline{Z}})+E_2$ and $K_{\tilde{Z}}+E_3=\sigma^*K_{Z'}+E_4$, where $E_i$ are $\mathbb Q$-effective $\mu$-exceptional divisors for $1\leq i\leq 4$ . Thus $ \vol(\overline{Z}, K_{\overline{Z}}+\Delta_{\overline{Z}})=\vol(\tilde{Z}, \mu^*(K_{\overline{Z}}+\Delta_{\overline{Z}})+E_2+E_3)=\vol(\tilde{Z}, K_{\tilde{Z}}+\Delta_{\tilde{Z}}+E_1+E_3)\geq \vol(Z', K_{Z'})$.
\end{proof}

The reason that we are interested to know the lower bound of the restricted volume is that we like to apply Helmke's induction to cut down log canonical centers.

We briefly recall Helmke's work  \cite{Hel}. Let  $D$ be an ample effective $\mathbb Q$-divisor on a smooth projective variety $X$ of dimension $n$. We assume that $x\in D$ and $\mathrm{mult}_xD=m>0$. Let $c:=\mathrm{lct}(D, x)$ be the log canonical threshold of $D$ at $x$, namely $c$ is the maximal rational number such that $(X, cD)$ is log canonical at $x$. Let $Z$ be the minimal log canonical center of $(X, cD)$ through $x$ and denote $d=\dim Z$. Helmke and Ein independently introduced the local discrepancy\footnote{This is called the deficit of $(X, cD)$ at $x$ in \cite{E} and \cite{YZ}.} $b_x(X, cD)$ of $(X, cD)$ at $x$, which is   the rational number
\begin{eqnarray*}\mathrm{max}\{\mathrm{mult}_x E\mid &&E \mathrm{ \;is\; an\; effective\; \mathbb Q-divisor}\; \\&& \mathrm{such\; that\;} (X, cD+E) \mathrm{ \;is\; log\; canonical \;at}\; x\}.
\end{eqnarray*}
It is known that $0\leq b_x(X, cD)\leq n-cm$, $b_x(X, cD)\leq d$ and $b_x(X, cD)=0$ iff $Z=\{x\}$ (see  \cite{Hel} and \cite{E}).

Helmke's induction can be summarized as follows (see \cite[Proposition 3.2 and Theorem 4.3]{Hel}).
\begin{prop}\label{helmke} Under the above assumption, assume that  $c<1$.
 \begin{itemize}
  \item[(1)] If $(D^{\dim Z}\cdot Z)>(\frac{b_x(X, cD)}{1-c})^{\dim Z}\mathrm{mult}_xZ$,   there exists a rational number $0<c'<1-c$, an effective $\mathbb Q$-divisor $D'\sim_{\mathbb Q}c'D$, such that $(X, cD+D')$ is log canonical at $x$ and the minimal log canonical center $Z'$ of $(X, cD+D')$ through $x$ is a proper subset of $Z$.
  \item[(2)] If $m\geq n$, let $c_1=c+c'$ and $D_1=cD+D'$. Then $\frac{b_x(X, D_1)}{1-c_1}<\frac{b_x(X, cD)}{1-c}\leq n$.
  \item[(3)] $\mathrm{mult}_xZ\leq \binom{n-\lceil b_x(X, cD) \rceil}{n-d}$.
 \end{itemize}
\end{prop}
\begin{rema}\label{minimal-center}
Let $(A, L)$ be a polarized abelian variety. Assume that an effective $\mathbb Q$-divisor $D\sim_{\mathbb Q}tL$ such that $(X, D)$ is log canonical around the neutral element $o$ and $Z$ is the minimal log canonical center through $o$. Choose an integer $m\geq 1$ such that $\cI_Z\otimes L^{\otimes m}$ is globally generated. Then by a standard argument (see for instance \cite[Proposition 2.3]{Kaw}), we know that for a small perturbation $D'=(1-\epsilon)D+\eta H$ of $D$, where $0<\epsilon, \eta<<1$  and $H\in |\cI_Z\otimes L^{\otimes m}|$   general, $Z$ is a connected component of $\mathrm{Nklt}(A, D')$.
In particular, when $Z=\{o\}$, we have $\beta(L)\leq r'(L)\leq t$. When $Z$ is an abelian subvariety of $A$ and $\beta(L|_Z)\leq t$, we also have $\beta(L)\leq t$ by Ito's work (see Proposition \ref{ito-abelian}).
\end{rema}
\subsection{A generic vanishing approach}

We recall some results from \cite{J}.

 \begin{lemm}\label{divisor}
Let $(A, L)$ be a polarized abelian variety.
Assume that $D\sim_{\mathbb Q} tL$ is an effective  $\mathbb Q$-divisor such that there exists an effective divisor $H\preceq D$, then $\cI_{K}\langle tL\rangle$ is GV, where $K$ is the neutral component of the kernel of the morphism
\begin{eqnarray*}&&\varphi_H: A\rightarrow \Pic^0(A)\\
&& x\rightarrow t_x^*\cO_A(H)\otimes \cO_A(-H).\end{eqnarray*}
In particular, when $H$ is ample, $\cI_o\langle tL\rangle $ is GV and thus $\beta(L)\leq t$.
\end{lemm}

\begin{proof} This follows directly from the proof of (1) of \cite[Proposition 4.1]{J}.
\end{proof}

The following theorem is essentially the main result of \cite{J}.
\begin{theo}\label{jiang-theorem}
Let $(A, L)$ be a polarized abelian variety. Assume that
\begin{itemize}
\item there exists an irreducible normal subvariety $Z$ of $A$ such that $\cI_Z\langle t_0L\rangle$ is IT$^0$ for some positive rational number $t_0$ and a smooth model of $Z$ is of general type;
\item and there exists an effective $\mathbb Q$-Weil divisor $D_Z$ on $Z$ and an effective $\mathbb Q$-Cartier divisor $V_Z$ such that $K_Z+D_Z$ is $\mathbb Q$-Cartier and $t_0L|_Z\sim_{\mathbb Q}2(K_Z+D_Z)+V_Z$,
\end{itemize}
we have $\beta(L)\leq t_0$.

\end{theo}
Since   \cite[Theorem 1.4]{J} is not stated this way, let's briefly recall the proof for readers' convenience.
\begin{proof}It suffices to show that $\cI_0\langle t_0L\rangle$ is GV. After a translation, we may assume that $o\in Z$ is a smooth point of $Z$. We have the short exact sequence $$0\rightarrow \cI_Z\rightarrow \cI_o\rightarrow \cI_{o, Z}\rightarrow 0. $$ By the first condition, it suffices to  show that $\cI_{o, Z}\langle t_0L\rangle$ is GV.

We may also assume that $Z$ is smooth (one may check \cite[Subsection 5.1]{J} for the full argument). Then since $Z$ is of general type, $\omega_Z\otimes\cI_{o, Z}$ is GV (see \cite[Lemma 2.6]{J}). Then one take an integer $M$ sufficiently divisible such that $Mt_0\in \mathbb Z$ and consider the multiplication-by-$M$ map $\pi_M: A\rightarrow A$ and denote by $Z^M$ the inverse image $\pi_M^{-1}(Z)$. By the second condition, $\cO_{Z^M}(M^2t_0L-K_{Z^M})$ has a non-trivial section $s$ and via this section, we have a short exact sequence $$0\rightarrow K_{Z^M}\xrightarrow{\cdot s} \cO_{Z^M}(M^2t_0L)\rightarrow \cQ\rightarrow 0.$$ We then check that all terms in this short exact sequence is GV. We may assume that the zero locus of $s$ does not intersect with $o_M:=\pi_M^{-1}(o)$. Then we have another short exact sequence $$0\rightarrow \cI_{o_M}\otimes K_{Z^M}\xrightarrow{\cdot s} \cI_{o_M}\otimes\cO_{Z^M}(M^2t_0L)\rightarrow \cQ\rightarrow 0.$$ Note that $\cI_{o_M}\otimes K_{Z^M}=\pi_M^{*}(\cI_o\otimes K_Z)$ is GV and hence $\cI_{o_M}\otimes\cO_{Z^M}(M^2t_0L)$ is also GV. This implies that $\cI_{o,Z}\langle t_0L\rangle$ is GV.
\end{proof}

We shall apply Theorem \ref{jiang-theorem} in the following way.
\begin{prop}\label{1/2} Let $(A, L)$ be a polarized simple abelian variety. Assume that an effective $\mathbb Q$-divisor $D\sim_{\mathbb Q}t_0L$ and $\mathrm{lct}(D, o)<\frac{1}{2}$. Then $\beta(L)\leq 2\mathrm{lct}(D, o)t_0<t_0$.
\end{prop}
\begin{proof}
We take $c=\mathrm{lct}(D)$ the global log canonical threshold of $D$, i.e. $c$ is the maximal rational number such that $(A, cD)$ is a log canonical pair. Then $c\leq \mathrm{lct}(D, o)<\frac{1}{2}$. After a small perturbation of $cD$ and taking a translation, we may assume that the log canonical pair $(A, cD)$ has only one log canonical center $Z$ which is thus an irreducible normal subvariety of $A$ and $o\in Z$.   Thus the multiplier ideal sheaf $\mathcal{J}(A, cD)=\cI_Z$ and by Nadel's vanishing, we have $\cI_Z\langle   2ct_0L\rangle$ is IT$^0$. Since $A$ is simple, a smooth model of $Z$ is of general type.

On the other hand, by the main theorem of \cite{FG}, we know that there exists an effective $\mathbb Q$-Weil divisor $D_Z$ on $Z$ such that $K_Z+D_Z$ is $\mathbb Q$-Cartier and $cD|_Z\sim_{\mathbb Q}K_Z+D_Z$. Thus by Theorem \ref{jiang-theorem}, $ \beta(L)\leq 2ct_0<t_0$.
\end{proof}

\subsection{Canonical volumes of subvarieties of abelian varieties}
Barja, Pardini, Stoppino studied higher dimensional Severi inequalities for varieties of maximal Albanese dimension in \cite[Corollary 5.6]{BPS}. One of their main result is the following.
\begin{theo}\label{BPS}Let $a: X\rightarrow A$ be a morphism from a smooth projective variety of general type of dimension $d\geq 2$ to an abelian variety $A$. Assume that $a: S\rightarrow a(S)$ is of degree $1$. Then $\vol(K_X)\geq \frac{5}{2}d!\chi(\omega_X) .$
\end{theo}

On the other hand, Pareschi and Popa generalized the Castelnuovo-De Franchis inequality in \cite{PP1}.
\begin{theo}\label{cd-inequality}Let $a: X\rightarrow A$ be a morphism from a smooth projective variety of general type of dimension $d\geq 2$ to an abelian variety $A$ of dimension $g$. Assume that  $A$ is simple and $a$ is generically finite from $X$ onto its image, $\chi(\omega_X)\geq g-d$.
\end{theo}

\begin{rema}\label{LP}
In \cite{LP}, Lazarsfeld and Popa conjectured that under that assumption of Theorem \ref{cd-inequality}, $\chi(\omega_X)> g-d$ when $g$ is large compared to $\chi(\omega_X)$. They verified that $\chi(\omega_X)\geq 3$ when $g\geq 5$ and $g-d=2$ (\cite[Proposition 4.10]{LP})
\end{rema}
Combining Theorem \ref{BPS} and \ref{cd-inequality}, we have an estimation of the canonical volume of irreducible subvarieties of simple abelian varieties.
\begin{coro}
Let $Z$ be an irreducible subvariety of dimension $d$ of a simple abelian variety $A$ of dimension $g$. Let $\rho: X\rightarrow Z$ be a desingularization. Then $\vol(K_X)\geq \frac{5}{2}d!(g-d)$.
\end{coro}

  When $X$ is a surface, we can improve slightly the above inequality.
\begin{prop}\label{variant}Let $a: S\rightarrow A$ be  a morphism from a smooth projective surface of general type to an abelian variety $A$ of dimension $g\geq 3$. Assume that $a: S\rightarrow a(S)$ is of degree $1$. Then $\vol(K_S)\geq 5\chi(\omega_S)+1$.
\end{prop}

\begin{proof}
We follows the proof of \cite[Corollary 5.6]{BPS}, where it was proved that
$\vol(K_S)\geq 5\chi(\omega_S)$. It suffices to show that the equality can never hold.

We may assume that $S$ is minimal and $a^*: \Pic^0(A)\rightarrow \Pic^0(S)$ is injective.

 Fix an ample divisor $H$ on $A$. Considering the continuous rank functions: $F(t):=h^0_a(A, a_*\omega_S\otimes H^t)$ and $G(t):=h^0_a(A, a_*(\omega_S^{\otimes 2})\otimes H^{2t})$. Following the proof in \cite[Theorem 5.5]{BPS}, we know that $D^{-1}G(t)\geq 6 D^{-1}F(t)$ for $t\leq 0$, where $D^{-1}G(t)$ and $D^{-1}F(t)$ are respectively the left derivative of the functions $F$ and $G$ at $t$. We also observe that $G(t)=F(t)=0$ when $t<<0$.
 Hence if
$G(0)=6F(0)$, i.e. $\vol(K_S)= 5\chi(\omega_S)$, $G(t)=6F(t)$ for all $t\leq 0$.

Since $a_*(\omega_S^2)$ is an IT$^0$ sheaf on $A$,
for $-\epsilon <t<0$, $a_*(\omega_S^2)\otimes H^{2t}$ remains to be IT$^0$ (see \cite[Theorem 5.2]{JP}). Thus $$G(t)=\chi(S, 2K_S+2tH)=2(H^2)_St^2+3(K_S\cdot H)_St+K_S^2+\chi(\omega_S)$$ is a degree $2 $ polynomial function for $-\epsilon <t<0.$

On the other hand, \begin{eqnarray*}F(t)&=&\chi(S, K_S+tH)+h^1_a(A, a_*\omega_S\otimes H^t)-h^2_a(A, a_*\omega_S\otimes H^t)\\
&=&\frac{1}{2}(H^2)_St^2+\frac{1}{2}(K_S\cdot H)_St+\chi(\omega_S)_S+h^1_a(A, a_*\omega_S\otimes H^t)\\&& - h^2_a(A, a_*\omega_S\otimes H^t).\end{eqnarray*} We recall the formula in \cite{JP} here:
\begin{eqnarray*}h^i_a(A, a_*\omega_S\otimes H^t)=\frac{(-t)^g}{\chi(H)}\chi(\phi_H^*R^i\Phi_{\cP}(a_*\omega_S)\otimes H^{-\frac{1}{t}}),
\end{eqnarray*}for $-\epsilon <t<0.$

Since $a^*: \Pic^0(A)\rightarrow \Pic^0(S)$ is injective, we know that $R^2\Phi_{\cP}(a_*\omega_S)$ is the skyscraper sheaf at $0_{\hat{A}}$. Hence $h^2_a(A, a_*\omega_S\otimes H^t)=\chi(H)(-t)^g$ for $-\epsilon <t<0.$

 Similarly, since $S$ is of general type, $a_*\omega_S$ is a M-regular sheaf on $A$ (see Subsection 2.3). Then $\codim_{\Pic^0(A)}V^1(a_*\omega_S)\geq 2$. We also know that the support of $R^1\Phi_{\cP}(a_*\omega_S)$ is contained in $V^1(a_*\omega_S)$.
Consider the Chern characters of $R^1\Phi_{\cP}(a_*\omega_S)$, we have $\mathrm{ch}_i(\varphi_H^*R^1\Phi_{\cP}(a_*\omega_S)=0\in H^{2i}(A, \mathbb Q)$
 for $i=0, 1$. If the codimension of the support of $R^1\Phi_{\cP}(a_*\omega_S)$ is equal to $2$, we may assume that $\cZ_1,\ldots, \cZ_s$ are the codimension-$2$  components of its  support and  $a_i>0$ is the rank of $R^1\Phi_{\cP}(a_*\omega_S)$ at the generic point of $\cZ_i$. Then
$\mathrm{ch}_2(\varphi_H^*R^1\Phi_{\cP}(a_*\omega_S)=\sum_ia_i[\varphi_H^{-1}\cZ_i]$.
 Hence  $$h^1_a(A, a_*\omega_S\otimes H^t)=\alpha_{2}t^2+\sum_{3\leq i\leq g}\alpha_it^i$$ for $-\epsilon <t<0,$ where $\alpha_2\geq 0$.

We compare the coefficient of $t^2$ for $F(t)$ and $G(t)$, which are respectively $\frac{1}{2}(H^2)_S+\alpha_2$ and $2(H^2)_S$ and conclude that $G(t)\neq 6F(t)$.
\end{proof}

 The following corollary seems to be new and it would be interesting to characterize the surfaces where the equality holds.
\begin{coro}\label{volume}
Let $Z$ be a smooth projective surface of general type. Assume that  $q(Z)\geq 4$ and $Z$ is of maximal Albanese dimension, either $Z$ is birational to a product of two smooth projective curve of genus $2$ or $\vol(Z, K_Z)\geq 16$.
\end{coro}

\begin{proof}
 If $\chi(\omega_Z)\geq 3$, we apply Proposition \ref{variant} to conclude that $$\vol(Z, K_{Z})\geq 16.$$ We then consider the case that $\chi(\omega_Z)=1$ or $2$.

When $\chi(\omega_Z)=1$, by the main result of \cite{HP} and the assumption that $q(Z)\geq 4$, we know that $q(Z)=4$ and $Z$ is birational to a product of two smooth projective curve of genus $2$.

 We now consider the case $\chi(\omega_{Z})= 2$ and $q(Z)\geq 4$.   When $q(Z)\geq 5$, then $p_g(Z)=q(Z)+1<2q(Z)-3$.
By the   Castelnuovo-de Francis theorem (see for instance \cite[Corollary 4.1]{PP1}), there exists a fibration $f: Z\rightarrow C$ to a smooth projective curve of genus $\geq 2$. We denote by $F$ a general fiber of $f$. Then $g(C)+g(F)\geq q(Z)\geq 5$. Take $Q\in \Pic^0(Z)$ a general torsion line bundle, we know by generic vanishing   that $h^0(K_{Z}\otimes Q)=\chi(\omega_{Z})=2$. On the other hand, by Viehweg's weak positivity, we know that $f_*(\omega_{Z/C}\otimes Q)$ is a nef vector bundle of rank $ g(F)-1$. By Riemann-Roch, we conclude that $h^0(K_{Z}\otimes Q)\geq (g(F)-1)(g(C)-1)$. The only possibility is that  $\{g(C), g(F)\}=\{2, 3\}$. Then $Z$ is biratioan to $C\times F$ by \cite[the Lemme in the appendix]{D2} and hence $\vol(Z, K_Z)=16$.
When  $q(Z)=4$, then $p_g(Z)=5$. By the main result of \cite{BNP},  we have  $16\leq \vol(Z, K_Z)\leq 18$.
\end{proof}

\subsection{Intersections with abelian subvarieties}
When we apply Ito's approach to study basepoint freeness thresholds of polarized abelian varieties $(A, L)$,  it may happen that the log canonical center $Z$ of a pair $(A, D)$ is an abelian subvariety of $A$ or is a subvariety fibred by an abelian subvariety of $A$.

Ito realized that one can deal with these cases with Poincar\'e's reducibility theorem  for polarized abelian varieties. When $Z$ is an abelian subvariety, we have the following result \cite[Proposition 6.6]{Ito}.
\begin{prop}\label{ito-abelian}
Let $(A, L) $ be a polarized abelian variety. Assume that there exists an effective $\mathbb Q$-divisor $D$ such that an abelian subvariety $B$ is an irreducible component of $\mathrm{Nklt}(A, D)$ and $tL-D$ is an ample $\mathbb Q$-divisor. Assume furthermore that $\beta(L|_B)<t$. Then $\beta(L)<t$.
\end{prop}


When $Z$ is fibred by an abelian subvariety $B$, it is important to estimate the intersection number $(L^{\dim B}\cdot B)$. I learned the following result from a private communication of Atsushi Ito.

\begin{lemm}\label{quotient}Let $D$ be an ample  $\mathbb Q$-divisor on $A$ and let $B$ be a abelian subvariety of $A$ of dimension $d$. Let $\varphi: A\rightarrow A/B$ be the quotient morphism. Then there exists an  ample $\mathbb Q$-divisor $H_B$ on $A/B$ such that $D-\varphi^*H_B$ is an effective nef $\mathbb Q$-divisor and $D^g=\binom{g}{d}(D^d\cdot B)(H_B^{g-d})_{A/B}$.
\end{lemm}

\begin{proof}Let $M$ be a positive sufficiently divisible integer such that $MD$ is an integral divisor and $L=\cO_A(MD)$ is an ample line bundle on $A$. By Poincar\'e's reducibility, there exists an abelian subvariety $K$ of $A$ such that the addition map $\mu: K\times B\rightarrow A$ is an isogeny and $\mu^*L\simeq L|_K\boxtimes L|_B$. Then one consider the natural isogeny $\mu_K: K\rightarrow A/B$. Note that $L|_K$ cannot descend to $A/B$ in general but there exists an effective $\mathbb Q$-divisor $D_B$ on $A/B$ such that $L|_K$ is algebraically equivalent to $\mu_K^*D_B$ as $\mathbb Q$-divisors. Let $H_B=\frac{1}{M}D_B$ and it is easy to check that $H_B$ satisfies the desired properties.\end{proof}
\begin{coro}\label{bigger}Let $(A, L)$ be a polarized abelian variety of dimension $g$ with $n(L)<\frac{1}{p+2}$. Assume that  Conjecture \ref{P2} holds in dimension $\leq g-1$, then either \begin{eqnarray}\label{newbound}(L^d\cdot B)\geq \frac{(L^g)}{\binom{g}{d}\big((p+2)(g-d)\big)^{g-d}} \end{eqnarray} for all abelian subvariety $B$ of dimension $1\leq d\leq g-1$ or $\beta(L)<\frac{1}{p+2}$.
\end{coro}
\begin{rema}Let $0<d<g$ be positive integers. Then $$\frac{g^g}{d^d\binom{g}{d}}>(g-d)^{g-d}.$$ More precisely, by Stirling's formula, for any positive integer $n$, we have $$\sqrt{2\pi}n^{n+\frac{1}{2}}e^{-n+\frac{1}{12n+1}}<n!<\sqrt{2\pi}n^{n+\frac{1}{2}}e^{-n+\frac{1}{12n}}.$$ Thus $$\frac{g^g}{d^d\binom{g}{d}}>(g-d)^{g-d}\sqrt{\frac{2\pi d(g-d)}{g}}.$$
\end{rema}
\begin{proof}
If the inequalities (\ref{newbound}) fails for some abelian subvarieties, we pick $B$ such that $\dim B=d$ is maximal and $(L^d\cdot B)< \frac{(L^g)}{\binom{g}{d}\big((p+2)(g-d)\big)^{g-d}}$. Let $D=\frac{1}{p+2}L$ and apply Lemma \ref{quotient}. There exists an ample $\mathbb Q$-divisor $H_{B}$  such that $D-\varphi^*H_{A/B}$ is a nef $\mathbb Q$-divisor and $$(H_{B}^{g-d})_{A/B}=\frac{(D^g)}{\binom{g}{d}(D^d\cdot B)}=\frac{(L^g)}{\binom{g}{d}(p+2)^{g-d}(L^d\cdot B)}>(g-d)^{g-d}.$$ Moreover, by the argument of Lemma \ref{quotient}, $(H_{B}^r\cdot K)_{A/B}=\frac{(D^{d+r}\cdot \varphi^{-1}(K))}{\binom{d+r}{d}(D^d\cdot B)}$.
 By the maximality of $B$, we conclude that for any abelian subvariety $K$ of $A/B$ whose dimension is $r>0$, we have
\begin{eqnarray*}(H_{B}^r\cdot K)_{A/B}=\frac{(D^{d+r}\cdot \varphi^{-1}(K))}{\binom{d+r}{d}(D^d\cdot B)}>&&\frac{ \binom{g}{d}(g-d)^{g-d}}{\binom{d+r}{d}\binom{g}{d+r}(g-d-r)^{g-d-r}}\\
&&=\frac{ (g-d)^{g-d}}{(g-d-r)^{g-d-r}\binom{g-d}{r}}\\&&>r^r.
\end{eqnarray*}
Since Conjecture \ref{P2} holds in dimension $g-d$, $\cI_{o_{A/B}}\langle (1-\epsilon)H_B\rangle$ is GV for some $0<\epsilon<<1$. Thus $\cI_B\langle (1-\epsilon)\varphi^*H_B\rangle$ is also GV. Since $D-\varphi^*H_B$ is an effective nef $\mathbb Q$-divisor, $\cI_B\langle (1-\epsilon)D\rangle$ is GV. We then consider the short exact sequence $$0\rightarrow \cI_B\rightarrow \cI_o\rightarrow \cI_{o, B}\rightarrow 0. $$  Since by the assumption that Conjecture \ref{P2} holds in dimension $d$, $\cI_{o, B}\langle (1-\epsilon)D\rangle$ is also GV and hence $\cI_o\langle (1-\epsilon)D\rangle$ is GV. Thus $\beta(L)<\frac{1}{p+2}$.
\end{proof}

\begin{lemm}\label{product-volume}Let $D$ be an effective $\mathbb Q$-divisor on an abelian variety $A$. Assume that $(A, D)$ is log canonical at $o$ and  the minimal log canonical center $Z$  through $o$ is a subvariety of dimension $d$ fibred by an abelian subvariety $B$ of dimension $d'$ such that the desingularization of $Z/B$ is of general type. Then $$(D^d\cdot Z)\geq \binom{d}{d'}(D^{d'}\cdot B)\vol(K_{(Z/B)'}),$$ where $(Z/B)'$ is a smooth model of the quotient $Z/B$.
\end{lemm}
\begin{proof} We consider the quotient morphism $A\rightarrow A/B$. We may assume that $MD$ is an integral divisor corresponding to a line bundle $L$ for some integer $M>0$. By Poincar\'e's reducibility theorem (see for instance \cite[Corollary 5.3.6]{BL}), there exists an abelian subvariety $K$ of $A$ such that the natural morphism $K\rightarrow A/B$ is an isogeny and the addition morphism $\pi: K\times B\rightarrow A$ induces an isogeny of polarized abelian varieties $(K, L_K)\times (B, L_B)\rightarrow (A, L)$ where $L_K$ and $L_B$ are respectively the restriction of $L$ on $K$ and $B$.

  Note that $\tilde{Z}:=Z\times_A(K\times B)$ is isomorphic to the product $\tilde{Z/B}\times B$, where $\tilde{Z/B}=(Z/B)\times_{A/B}K$. Thus,
\begin{eqnarray*}
&&(D^d\cdot Z)=\frac{1}{M^d}(L^d\cdot Z)=\frac{1}{M^d\deg\pi}((L_K\boxtimes L_B)^d\cdot\tilde{Z})\\
&=&\frac{1}{M^d\deg\pi}\big((L_K\boxtimes L_B)^d\cdot(\tilde{Z/B}\times B)\big)\\&=&\frac{\binom{d}{d'}}{M^{d-d'}\deg\pi}(D^{d'}\cdot B)(L_K^{d-d'}\cdot \tilde{Z/B}).
\end{eqnarray*}
Let $\nu: \overline{Z}\rightarrow Z$ be its normalization. By Proposition \ref{Jx}, $\nu^*D\sim_{\mathbb Q}K_{\overline{Z}}+V_{\overline{Z}}$ for some effective $\mathbb Q$-divisor $V_{\overline{Z}}$. Consider the pull-back of this $\mathbb Q$-linear equivalence on the normalization $ \overline{\tilde{Z}}$ of $\tilde{Z}$ and then restricting it to a general fiber of $ \overline{\tilde{Z}}\rightarrow B$, we see that $\varsigma^*(\frac{1}{M}L_{K})\sim_{\mathbb Q}K_{\widehat{Z/B}}+V_{\widehat{Z/B}} $, where $\varsigma:\widehat{Z/B}\rightarrow \tilde{Z/B}$ is the normalization and $V_{\widehat{Z/B}}$ is an effective $\mathbb Q$-divisor on $\widehat{Z/B}$. Thus $(L_K^{d-d'}\cdot \tilde{Z/B})= M^{d-d'}\vol(\widehat{Z/B}, K_{\widehat{Z/B}}+V_{\widehat{Z/B}})\geq (\deg\pi )M^{d-d'}\vol(K_{(Z/B)'})$.
\end{proof}
\section{Very general polarized abelian varieties}
In \cite{J}, we applied Helmke's induction to confirm that Conjecture \ref{P2} holds for Hodge theoretically very general polarized abelian varieties with special polarizations.
We now show that the calculation indeed implies that Conjecture \ref{P2} holds for almost  all generic polarized abelian varieties in fixed dimensions.

Recall that we say a polarized abelian variety $(A, L)$ is Hodge theoretically very general if $\dim_{\mathbb Q}H^{k,k}(A, \mathbb Q)=1$ for all $1\leq k\leq g-1$. We observe that by Hard Lefschetz and Poincar\'e duality,
$(A, L)$ is Hodge theoretically very general iff $\dim_{\mathbb Q}H^{\lfloor \frac{g}{2}\rfloor,\lfloor \frac{g}{2}\rfloor}(A, \mathbb Q)=1$. We also observe that when $(A, L)$ is Hodge theoretically very general, $A$ is a simple abelian variety and hence $n(L)=\frac{g}{\sqrt[g]{(L^g)}}$. Note that in order to compare $\beta(L)$ and $n(L)$, we can assume that $L$ is primitive.

\begin{theo}\label{finitely-many}Let $(A, L)$ be a Hodge theoretically very general polarized abelian variety of type $l=(1,\delta_2,\ldots,\delta_g)$. Assume that \begin{eqnarray*}\label{bound}\delta:=\delta_2\cdots\delta_g\geq \mathrm{ max}\big\{\frac{(k(k+1)\cdots(g-1))^{\frac{g(g-1)}{k}}}{(g!)^{\frac{(g-1)(g-k)}{k}}}\mid 2\leq k\leq g-2\big\},\end{eqnarray*}
$$\beta(L) \leq n(L).$$
\end{theo}
\begin{proof}
We apply Ito's strategy that it suffices to show $r'(L)\leq n(L)$. Let $t\in (n(L), n(L)+\epsilon)$ be a rational number, where $0<\epsilon<<1$ and denote by $D\sim_{\mathbb Q} tL$ an effective rational number such that $\mathrm{mult}_0D>g$ (there exists such $D$ since $((tL)^g)>g^g$). We then just need to apply Helmke's induction \ref{helmke} to get a divisor $D'\sim_{\mathbb Q}cD$ with $0<c<1$ and the neutral element $o$ of $A$ is a minimal log canonical center of $(A, D')$. By Proposition \ref{helmke}, it suffices to verify that $(D^{k}\cdot Z)>\binom{g-1}{g-k}g^k$ for any irreducible subvariety $Z$ of dimension $k$ for $1\leq k\leq g-1$.

Recall that $\beta(L)\leq 1$ and equality holds if and only if $|L|$ is basepoint free. Thus we may assume that $(L^g)>g^g$, i.e. $n(L)<1$. Then since $[L]$ is the generator of $\mathrm{NS}(A)=\Pic(A)/\Pic^0(A)$, we can assume that $Z$ is of codimension $\geq 2$.  Since $A$ is simple, by the main result of \cite{Deb},  for any irreducible curve $Z$ of $A$, $(L\cdot Z)>\sqrt[g]{(L^g)}>g$. Thus it suffices to verify that $(D^{k}\cdot Z)>\binom{g-1}{g-k}g^k$ for any irreducible subvariety $Z$ of dimension $k$, where $2\leq k\leq g-2$.

For any irreducible subvariety $Z$ of codimension $k$, we denote by $[Z]\in H^{k, k}(A, \mathbb Z)\subset H^{2k}(A, \mathbb C)$ its cohomology class. Then $[Z]$ is a positive integral multiple of $\frac{L^{g-k}}{(g-k)!\delta_2\cdots\delta_{g-k}}$ by the assumption that $A$ is Hodge theoretically very general. Thus
\begin{eqnarray*}(D^k\cdot Z)> \frac{g^k}{\sqrt[\frac{k}{g}]{(L^g)}}(L^k\cdot Z) \geq g^k\frac{\sqrt[\frac{g-k}{g}]{(L^g)}}{(g-k)!\delta_2\cdots\delta_{g-k}}.
\end{eqnarray*}
We just need to verify that \begin{eqnarray}\label{hodge-computation} \frac{\sqrt[\frac{g-k}{g}]{(L^g)}}{(g-k)!\delta_2\cdots\delta_{g-k}}\geq \binom{g-1}{g-k}.\end{eqnarray} Note that $\delta_2\cdots\delta_{g-k}\leq \delta^{\frac{g-k-1}{g-1}}$. Thus \begin{eqnarray*} \frac{\sqrt[\frac{g-k}{g}]{(L^g)}}{(g-k)!\delta_2\cdots\delta_{g-k}} \geq \delta^{\frac{g-k}{g}-\frac{g-k-1}{g-1}}\frac{\sqrt[\frac{g-k}{g}]{g!}}{(g-k)!}\\=\delta^{\frac{k}{g(g-1)}}\frac{\sqrt[\frac{g-k}{g}]{g!}}{(g-k)!}.
\end{eqnarray*}
The assumption on $\delta$ makes sure that $  \frac{\sqrt[\frac{g-k}{g}]{(L^g)}}{(g-k)!\delta_2\cdots\delta_{g-k}}\geq \binom{g-1}{g-k}$ for $2\leq k\leq g-2$.
\end{proof}
By the above theorem, we see that in order to verify Conjecture \ref{P1} and Conjecture \ref{P2} for very general polarized abelian varieties, it suffices to check finitely many families in each dimension.
By the same computation, we finish the proof of Theorem \ref{thm1}.
\begin{proof}
When $g=4$, the assumption in Theorem (\ref{finitely-many}) is simply that $\delta=h^0(A, L)\geq \frac{6^6}{24^3}=\frac{27}{8}$. On the other hand, for $\delta\leq 3$, $n(L)=\frac{4}{\sqrt[4]{24\delta}}>1$ we know that $|L|$ has  basepoints and $\beta(L)=1$.

When $g=5$, we repeat the argument in the proof of Theorem (\ref{finitely-many}). We need to verify (\ref{hodge-computation}) for $g=5$ and $k=2$ or $3$, which are $(L^5)^{\frac{3}{5}}\geq 24\delta_2\delta_3$ and $(L^5)^{\frac{2}{5}}\geq 12\delta_2$. It is easy to verify that both inequalities hold when $\delta_5\geq 5$ or $\delta_2\geq 3$. When $\delta_2=1$,  $\delta\geq 5$ implies that both inequalities hold and if $\delta<5$, $\beta(L)=1$. When $\delta_2=2$ and $\delta_5\leq 4$, the above inequalities fail only when the polarization  is of type $(1,2,2,2,2)$. But in this case, $n(L)=\frac{5}{\sqrt[5]{120\delta(L)}}>1$, thus we still have $\beta(L)\leq 1<n(L)$.
\end{proof}

\begin{rema}
Let $(A, L)$ be a Hodge theoretically very general polarized abelian variety of dimension $6$. The same methods shows that $\beta(L)\leq n(L)$ unless the polarization type is proportional to $(1, 3, 3,3,3,3)$.
\end{rema}
 \section{Abelian fourfolds}

\subsection{The proof of Theorem \ref{thm2}}
\begin{proof}
It suffices to show that $\beta(L)<\frac{1}{p+2}$ or equivalently $\cI_o\langle \frac{1}{p+2}L\rangle$ is IT$^0$, where $o$ is the neutral element of $A$.

 Since $(L^4)>((2+\frac{4}{\sqrt{3}})(p+2))^4$, there exists an effective $\mathbb Q$-divisor $D\sim_{\mathbb Q}\frac{1}{p+2}L$ such that $\mathrm{mult}_o(D)=m>2+\frac{4}{\sqrt{3}}$. Let $c=\mathrm{lct}(D, o)\leq \frac{4}{m}$ be the log canonical threshold of $D$ at $o$ and let $Z$ be the minimal log canonical center of $(A, cD)$ through $o$. By Proposition \ref{1/2}, we may assume that $c>\frac{1}{2}$.

 {\bf Step 1.}
We first deal with the case that  $Z$ is a divisor. By Lemma \ref{divisor}, $\cI_K\langle \frac{c}{p+2}L\rangle$ is GV, where $K$ is the kernel of $\varphi_{Z_1}$ .
If $K$ is a point, we are done.
Otherwise, by the assumption that $(L^d\cdot B)>(d(p+2))^d$ for any abelian subvariety $B$ of dimension $1\leq d\leq 3$ and Ito's results \cite{Ito1, Ito}, $\cI_{o,K}\langle \frac{1}{p+2}L\rangle$ is IT$^0$. Thus from the short exact sequence $$0\rightarrow \cI_K\rightarrow \cI_o\rightarrow \cI_{o, K}\rightarrow 0,$$ we conclude that $\cI_o\langle \frac{1}{p+2}L\rangle$ is IT$^0$.

{\bf Step 2.}
If $Z$ is a curve, by Proposition \ref{helmke}, $Z$ is smooth at $o$ and as soon as $(D\cdot Z)>4\geq \frac{b_o(A, cD)}{1-c}$, there exists $D_1\sim_{\mathbb Q}c_1D$ with $c<c_1<1$ such that $(A, D_1)$ is log canonical at $o$ and $\{o\}$ is a minimal log canonical center. Then we have by Remark \ref{minimal-center} that $\beta(L)\leq r'(L)\leq \frac{c_1}{p+2}<\frac{1}{p+2}$.

  By Corollary \ref{intersection}, we know that $$(cD\cdot Z)\geq 2g(\overline{Z})-2, $$ where $\nu: \tilde{Z}\rightarrow Z$ is the normalization. Thus we are done once $g(\overline{Z})\geq 3$.

  When $g(\overline{Z})=1$, $Z=\overline{Z_1}$ is an elliptic curve. We conclude again by Proposition \ref{ito-abelian}, Remark \ref{minimal-center} and the assumption that $(L\cdot Z_1)>p+2$ and thus $\beta(L|_{Z_1})<\frac{1}{p+2}$.

  If $g(\overline{Z})=2$, $Z$ generates an abelian surface $B$ of $A$. By Corollary \ref{bigger}, we may assume that $(D^2\cdot B)\geq \frac{4^4}{\binom{4}{2}2^2}= \frac{32}{3}$.  Then by Hodge index theorem, $$(D\cdot Z)_B\geq \sqrt{(Z^2)_B(D^2\cdot B)} >8/\sqrt{3}>4.$$

{\bf Step 3.}
When $Z$ is a surface, we need to apply Helmke's induction. We may assume that $Z$ is not an abelian surface, otherwise we conclude directly by Proposition \ref{ito-abelian} and Ito's work in \cite{Ito1}. We know that $\mathrm{mult}_oZ\leq 3$ by Proposition \ref{helmke}. Let $\mu: \tilde{Z}\rightarrow \overline{Z}\rightarrow Z$ be the minimal resolution of the normalization $\overline{Z}$ of $Z$. Since $Z$ is not an abelian variety, $\tilde{Z}$ is a surface of maximal Albanese dimension of Kodaira dimension $\geq 1$. Since there exists no rational curves on $\overline{Z}$, $\tilde{Z}$ is minimal.

\textbf{Claim:}  $((cD)^2\cdot Z)\geq 16$, when $\mathrm{mult}_oZ=3$; $((cD)^2\cdot Z)\geq 16$ or $(D^2\cdot Z)>16\sqrt[3]{6}$, when $\mathrm{mult}_oZ=2$; $(D^2\cdot Z)>16$ when $Z$ is smooth at $o$.

When $\tilde{Z}$ is not of general type, $Z$  is fibred by an elliptic curve $E$ and since $Z$ has rational singularities around $o$, $Z$ is indeed smooth at $o$. We need to show that $(D^2\cdot Z)>16$.
Let $C=Z/E$ be the quotient and $\tilde{C}$  be the normalization of $C$. By Corollary \ref{bigger}, we may assume $(D\cdot E)\geq \frac{64}{27}$ and by Lemma \ref{product-volume}, we have $$(D^2\cdot Z)\geq 2(2g(\tilde{C})-2)(D\cdot E).$$  Thus we are done when $g(\tilde{C})\geq 3$. When $g(\tilde{C})=2$, $Z$ generates an abelian $3$-fold $B$ of $A$ and $C\hookrightarrow B/E$ is an ample divisor. Thus $(C^2)_{B/E}=a\geq 2$ and hence $Z^2$ is algebraically equivalent to $aE$ as $1$-cycles of $B$.
By Corollary \ref{bigger}, we may assume that $(D^3\cdot B)\geq 4^3$. Then by Hodge index, \begin{eqnarray*}(D^2\cdot Z)=((D|_B)^2\cdot Z)_B&&\geq \sqrt{(D^3)_B(D|_B\cdot Z^2)_B}\\&& \geq \sqrt{4^3\cdot a\cdot \frac{64}{27}}\geq 16\sqrt{\frac{32}{27}}>16.\end{eqnarray*}

We now assume that $\tilde{Z}$ is of general type.

 If $q(\tilde{Z})\geq 4$. By Corollary \ref{volume}, $\vol(K_{\tilde{Z}})\geq 16$ or $\tilde{Z}\simeq C_1\times C_2$, where $C_i$ is a smooth projective curve of genus $2$. In the latter case $\rho$ is the normalization of $Z$ and since $Z$ is normal at $o$, it is smooth at $o$.
We apply Proposition \ref{Jx} and conclude that \begin{eqnarray*}(D^2\cdot Z)=(\rho^*(D)^2)_{\tilde{Z}}&&\geq (K_{\tilde{Z}}\cdot \rho^*D)_{\tilde{Z}} \\&&=2(C_1\cdot \rho^*D)_{\tilde{Z}} +2(C_2\cdot \rho^*D)_{\tilde{Z}} .\end{eqnarray*}
Note that the image of $C_i$ generates an abelian surface $B_i$ of $A$. By Corollary  \ref{bigger}, we may assume that $(D^2\cdot B_i)\geq \frac{32}{3}$. Thus by Hodge index, $(C_i\cdot \rho^*D)_{\tilde{Z}}\geq 8/\sqrt{3}$. Thus we have $(D^2\cdot Z)\geq 32/\sqrt{3}>16$.
In the former case, $ ((cD)^2\cdot Z)=\vol(Z, cD)\geq \vol(K_{\tilde{Z}})\geq 16$ by Corollary \ref{intersection}.

 If $q(\tilde{Z})=3$, $Z$ generates an abelian 3-fold $B\subset A$. Then $Z$ is an ample divisor of $B$. Moreover, in this case the embedded dimension of $Z$ at $o$ is at most $3$, thus $\mathrm{mult}_o(Z)\leq 2$ by the well-known facts about isolated rational surface singularities (see for instance \cite{A}). As before, by Corollary \ref{bigger}, we may assume that $(D^3\cdot B)\geq 4^3$. Since $Z$ is an ample divisor of $B$, $(Z^3)_B\geq 6$. Thus $$(D^2\cdot Z)\geq \sqrt[3]{(D^3\cdot B)^2(Z^3)_B}\geq 16\sqrt[3]{6}.$$

{\bf Step 4.}

If $Z$ is smooth at $o$, we have shown that $(D^2\cdot Z)>16$, thus by Proposition \ref{helmke}, there exists an effective $\mathbb Q$-divisor $D_1\sim_{\mathbb Q} c_1D$ with $c_1<1$ such that $(A, D_1)$ is log canonical at $o$ whose  minimal log canonical center $Z_1$ through $o$ is a proper subset of $Z$ and $\frac{b_o(A, D_1)}{1-c_1}\leq \frac{b_o(A, cD)}{1-c}\leq 4$. We then finish the proof by going back to Step 2.

If $Z$ is singular at $o$, $\mathrm{mult}_oZ= 3$. We have already seen that $((cD)^2\cdot Z)\geq 16$. In order to apply Helmke's induction, we need to verify that $$(D^2\cdot Z)>3\big(\frac{b_o(A, cD)}{1-c}\big)^2.$$
Note that $b_o(A, cD)\leq 4-cm<4-(2+\frac{4}{\sqrt{3}})c$. It is elementary to verify that $$\frac{16}{c^2}\geq 3\big(\frac{4-(2+\frac{4}{\sqrt{3}})c}{1-c}\big)^2$$ always holds. We then finish the proof as before.

If $\mathrm{mult}_oZ=2$ and  $((cD)^2\cdot Z)\geq 16$, we conclude as the multiplicity $3$ case. If $\mathrm{mult}_oZ=2$ and  $(D^2\cdot Z)\geq 16\sqrt[3]{6}$, we need to verify that $$(D^2\cdot Z)>2\big(\frac{b_o(A, cD)}{1-c}\big)^2.$$ Since $c>\frac{1}{2}$, we have $\frac{b_o(A, cD)}{1-c}<4-\frac{(\frac{4}{\sqrt{3}}-2)c}{1-c}<6-\frac{4}{\sqrt{3}}$. We then check that $16\sqrt[3]{6}>2(6-\frac{4}{\sqrt{3}})^2$.
 \end{proof}

\subsection{The proof of Theorem \ref{thm3}}

We apply the same strategy as in the proof of Theorem \ref{thm2}.

Fix an effective $\mathbb Q$-divisor such that $D\sim_{\mathbb Q}\frac{1}{p+2}L$ such that $\mathrm{mult}_oD>8$. Let $c_1=\mathrm{lct}(D, o)<\frac{5}{8}$ be the log canonical threshold of $D$ at $o$ and let $Z_1$ be the minimal log canonical center of $(A, c_1D)$ at $o$.  By Proposition \ref{1/2}, we may assume that $c_1\geq \frac{1}{2}$. Thus
  $$\frac{b_o(A, c_1D)}{1-c_1}\leq \frac{5-c_1\mathrm{mult}_oD}{1-c_1}<5-3\frac{c_1}{1-c_1}\leq 2.$$

If $Z_1$ is a divisor, we conclude by Lemma \ref{divisor}.

If $Z_1$ is a threefold, we apply Helmke's induction. Let $\rho: \tilde{Z_1}\rightarrow Z_1$ be a desingularization. Then $\vol(K_{\tilde{Z_1}})\geq \frac{5}{2}3!\times 3=45$ by Theorem \ref{BPS}, Theorem \ref{cd-inequality}, and Remark \ref{LP}. We also note that $\mathrm{mult}_oZ_1\leq 6$.
Note that $$(D^3\cdot Z_1)\geq \vol(K_{\tilde{Z_1}})/c_1^3>\frac{45}{(\frac{5}{8})^3}\geq 48>\big(\frac{b_o(A, c_1D)}{1-c_1}\big)^3\mathrm{mult}_{o}(Z_1).$$
Thus there exists an effective $\mathbb Q$-divisor  $ D_2\sim_{\mathbb Q}c_2D$ such that $c_1<c_2<1$ such that $(A, D_2)$ is log canonical at $o$, whose minimal lc center through $o$ is a proper subvariety $Z_2$ contained in $Z_1$, and $\frac{b_o(A, D_2)}{1-c_2}<\frac{b_o(A, c_1D)}{1-c_1}<2$.

When $Z_2$ is a surface, we have $\mathrm{mult}_oZ_2\leq 4$. Let $\tilde{Z}_2$ be its smooth model. Since $A$ is simple, $\tilde{Z}_2$ is smooth and hence $\vol(K_{\tilde{Z_2}})\geq 16$ by Corollary \ref{volume}. We then have
$$(D^2\cdot Z_2)>\vol(K_{\tilde{Z_2}})\geq 16>\big(\frac{b_o(A, D_2)}{1-c_2}\big)^2\mathrm{mult}_oZ_2.$$ Thus by Helmke's induction, we may assume that $Z_2$ is a curve. In this case,  we verify easily that $(D\cdot Z_2)>\frac{b_o(A, D_2)}{1-c_2}$.  Therefore, there exists an effective $\mathbb Q$-divisor  $D_3\sim_{\mathbb Q}c_3L$ with $c_2<c_3<1$ such that $(A, D_3)$ is log canonical at $o$ and $o$ is a minimal log canonical center of $(A, D_3)$. We then finish the proof of Theorem \ref{thm3}.

\subsection{The proof of Theorem \ref{thm4}}

The proof of Theorem \ref{thm4} is identical to that of Theorem \ref{thm3}.
We first observe that $\alpha_g\leq \alpha_{g,2}=\sqrt{\frac{5(g-2)}{g-1}}<\sqrt{5}$.
 We then fix an effective $\mathbb Q$-divisor $D\sim\frac{1}{p+2}L$ such that $\mathrm{mult}_oD>2g-\alpha_g$ and thus let $c<\frac{g}{2g-\sqrt{5}}$ be the log canonical threshold of $(A, D)$ at $o$. We may assume that $c\geq \frac{1}{2}$ by Proposition \ref{1/2}. Then $$\frac{b_o(A, cD)}{1-c}\leq \frac{g-c (\mathrm{mult}_oD)}{1-c}\leq g-(\mathrm{mult}_oD-g)<\alpha_g.$$ Let $Z$ be the minimal log canonical threshold of $(A, cD)$ through $o$. Let $d=\dim Z$. Then $\mathrm{mult}_oD\leq \binom{g-1}{g-d}$.

 If $d=g-1$, we conclude by Lemma \ref{divisor}. Thus we may assume that $1\leq d\leq g-2$.

 For a smooth model $\tilde{Z}$ of $Z$, by Theorem \ref{BPS}, Theorem \ref{cd-inequality}, Remark \ref{LP}, and Proposition \ref{variant}, we have  $\vol(K_{\tilde{Z}})\geq \frac{15}{2}(g-2)!$ when $d=g-2$, $\vol(K_{\tilde{Z}})\geq \frac{5}{2}d!(g-d)$ when $3\leq d\leq g-3$, $\vol(K_{\tilde{Z}})\geq 5(g-2)+1$ when $d=2$,  and $\vol(K_{\tilde{Z}})\geq 2(g-1)$. By Corollary \ref{intersection}, $(D^d\cdot Z)\geq \frac{1}{c^d}\vol(K_{\tilde{Z}})>\vol(K_{\tilde{Z}})$. We then have $$(D^d\cdot Z)> (\alpha_g)^d\mathrm{mult}_oD.$$

 We then repeatedly apply Helmke's induction and the above calculation to  cut down the log canonical centers and finish the proof.

\end{document}